\documentclass{amsart}[12pt,a4paper]

\usepackage[colorlinks=true,          
            linkcolor=blue,
            citecolor=red,
            urlcolor=blue]{hyperref}
\usepackage{amssymb,amsfonts,amsmath,hyperref,mathrsfs,latexsym,stmaryrd,graphicx,manfnt}
\DeclareMathAlphabet{\mathpzc}{OT1}{pzc}{m}{it}

\DeclareMathAlphabet{\mathpzc}{OT1}{pzc}{m}{it}

\makeatletter
\renewcommand*\env@matrix[1][c]{\hskip -\arraycolsep
  \let\@ifnextchar\new@ifnextchar
  \array{*\c@MaxMatrixCols #1}}
\makeatother

\newcommand{\diag}{\textrm{diag}}

 \newtheorem{thm}{Theorem}[section]
  \theoremstyle{plain}\newtheorem*{thma}{Theorem A}
 \newtheorem{Lemma}{Lemma}[section]
 
 \newtheorem{Prop}{Proposition}[section]
  \theoremstyle{definition}\newtheorem{Ex}{Example}[section]

\newcommand {\CC}{\mathbb{C}}

\newcommand{\QQ}{\mathbb{Q}}
\newcommand {\RR}{\mathbb{R}}

\newcommand{\ZZ}{\mathbb{Z}}

\newcommand {\bal}{\boldsymbol{\alpha}}
\newcommand {\bbe}{\boldsymbol{\beta}}
\newcommand {\bga}{\boldsymbol{\gamma}}
\newcommand {\bla}{\boldsymbol {\lambda}}

\newcommand {\bdU}{\boldsymbol{U}}

\newcommand {\bg}{{\bf g}}

\newcommand {\bs}{{\bf s}}
\newcommand {\bt}{{\bf t}}

\newcommand {\bI}{{\bf I}}

\newcommand {\bU}{{\bf U}}

\newcommand{\cB}{\mathcal{B}}

\newcommand{\cF}{\mathcal{F}}

\newcommand{\cH}{\mathcal{H}}

\newcommand{\cJ}{\mathcal{J}}
\newcommand{\cK}{\mathcal{K}}
\newcommand{\cL}{\mathcal{L}}

\newcommand {\cO}{\mathcal{O}}

\newcommand{\cV}{\mathcal{V}}

\newcommand{\cX}{\mathcal{X}}

\newcommand{\scB}{\mathscr{B}}
\newcommand{\scC}{\mathscr{C}}

\newcommand{\fb}{\mathfrak{b}}

\newcommand{\fg}{\mathfrak{g}}

\newcommand{\fl}{\mathfrak{l}}

\newcommand{\fo}{\mathfrak{o}}

\newcommand{\fs}{\mathfrak{s}}

\newcommand{\ft}{\mathfrak{t}}

\newcommand{\sym}{\textrm{Sym}}

\newcommand{\mhom}{\textrm{Hom}}
\newcommand {\mend}{\textrm{End}}

\newcommand {\send}{\underline{End} }

\newcommand{\tot}{\textrm{tot }}

\newcommand{\ctimes}{\otimes_\CC}

\newcommand {\io}{\iota}

\newcommand{\Prym}{{\bf Prym}}

\newcommand{\hookr}{\hookrightarrow}

\include{biblio}

\usepackage{latexsym}
\usepackage{mathrsfs}
\input xy
\xyoption{all}

\title{ Seiberg--Witten differentials on the Hitchin base}
\author[U. Bruzzo and P. Dalakov]{Ugo Bruzzo$^{\P\S\sharp\star}$ and Peter Dalakov$^{\dag\ddag}$}
\address{\small $^\P$ SISSA (Scuola Internazionale Superiore di Studi Avanzati)\\
Via Bonomea 265, 34136 Trieste, Italia\\
$^\S$ Departamento de Matem\'atica, Universidade Federal da Para\'iba, Jo\~ao Pessoa, Brazil\\
$^\sharp$ INFN (Istituto Nazionale di Fisica Nucleare), Sezione di Trieste \\
$^\star$ IGAP (Institute for Geometry and Physics), Trieste, Italy\\
$^\dag$ American University in Bulgaria, 2700 Blagoevgrad, G.Izmirliev Sq.1, Bulgaria\\
$^\ddag$ Institute of Mathematics and Informatics, Bulgarian Academy of Sciences, Sofia, Bulgaria}

\date{Revised 7 November 2023}

\subjclass[2010]{14D20, 14D07,  14H70 }
\begin{document}

  \begin{abstract} In this note we describe explicitly,  in terms of Lie theory and 
	    cameral data,  the covariant (Gauss--Manin) derivative of the Seiberg--Witten differential   
	    defined on the  weight-one variation of Hodge structures that exists on a Zariski open subset of the base of the Hitchin fibration.
      \end{abstract}

\dedicatory{Dedicated to Tony Pantev on the occasion of his 60th birthday.}
\maketitle
\setcounter{tocdepth}{2}
\tableofcontents

\eject
      \section{Introduction}

	    The base of the Hitchin integrable system (\cite{hitchin_sb}) supports a family of cameral curves, and, as a consequence,
	    carries various Hodge-theoretic and differential-geometric structures (\cite{donagi_markman}, \cite{donagi_markman_cubic}, \cite{markman_sw}).
	    In particular, the Zariski open subset of the base,
	    corresponding to smooth cameral curves with generic ramification carries a weight-one variation of Hodge structures (VHS) with a Seiberg--Witten differential.
	    Our goal in this note is to describe the covariant (Gauss--Manin) derivative of the Seiberg--Witten differential explicitly in terms of Lie theory and 
	    cameral data.
	    
	    We recall now the main ingredients and constructions, starting with the Hodge-theoretic ones.

	      Let $\cB$ be a complex manifold. 
	    Recall that a polarised
    $\RR$-VHS of weight $w\in\ZZ$ on $\cB$ consists of 
    data $(\cV,\nabla,\cV_\RR,\cF^\bullet,S)$, where:  
	\begin{itemize}
		\item $\cV$ is a holomorphic vector bundle on $\cB$ 
		\item $\nabla:\cV\to \cV\otimes\Omega^1_\cB$ is a flat (holomorphic) connection, called   \emph{the Gauss--Manin connection}  
		\item $\cV_\RR\subseteq \cV$ is a real,  $\nabla$-flat subbundle, satisfying  $\cV=\cV_\RR\otimes \cO_\cB$, called \emph{real structure}
		\item $\cF^\bullet$ is  a decreasing  \emph{filtration} of $\cV=\cF^0$, \emph{the Hodge filtration}
		\item $S:\cV\otimes\cV\to \scC^\infty_\cB$ is a non-degenerate, $(-1)^w$-symmetric, $\nabla$-flat pairing, $\RR$-valued on $\cV_\RR$, called  \emph{polarisation}  
	\end{itemize}
	such that 
	\begin{enumerate}
		\item $\nabla(\cF^p)\subseteq \cF^{p-1}\otimes \Omega_{\cB}^1$  \emph{Griffiths transversality}
		\item $\cV= \cF^p\oplus \overline{\cF^{w+1-p}}$  \emph{Hodge structure},
	\end{enumerate}
	or, in terms of the Hodge bundles  $\cH^{p,w-p}:= \cF^p\cap\overline{\cF^{w-p}}$,
      \begin{enumerate}
      	\item[(2')] $\cV=\bigoplus_{p}\cH^{p,w-p}$
      	\item[(3)] $S(\cF^p,\cF^{w+1-p})=0$
      	\item[(4)] $i^{2p-w}S(v,\overline{v})>0$ for $v\in \Gamma(\cH^{p,w-p})$, $v\neq 0$.
      \end{enumerate}
      
      The notions of  polarised $\ZZ$-VHS or $\QQ$-VHS are introduced analogously, by replacing $\cV_\RR$ with  appropriate  locally constant
      sheaves $\cV_\ZZ$ or $\cV_\QQ$ of $\ZZ$- or $\QQ$-modules, respectively.

      The prototypical example is that of a geometric VHS, i.e., one arising  from a family of  compact K\"ahler (e.g., projective) manifolds.

		  By Griffiths Transversality, $\nabla$ induces an $\cO_B$-module homomorphism
	\[
		\cF^p/\cF^{p+1}\longrightarrow \cF^{p-1}/\cF^{p}\otimes \Omega^1_\cB
	\]
	and hence, taking a direct sum over the different $p$,  an $\cO_B$-module homomorphism
	\[
		\theta=[\nabla]: \bigoplus_p \cF^p/\cF^{p+1} \longrightarrow  \left(\bigoplus_p \cF^p/\cF^{p+1} \right)\otimes \Omega^1_\cB,
	\]
	which satisfies
	$\theta\wedge \theta=0$.

	The pair 
		$\left(E= \bigoplus_p \cF^p/\cF^{p+1},\theta  \right)$
	is an example of a \emph{Higgs bundle} on $\cB$. This example played an important r\^ole in  Carlos Simpson's 
	 study of Higgs bundles on higher-dimensional varieties  (\cite{hbls}, \cite{simpson_uniformisation}).

	    Consider  a polarised $\ZZ$-VHS $(\cV,\nabla, \cV_\ZZ, \cV^\bullet, S,\ldots )$ of weight $w=1$. An
	\emph{abstract Seiberg--Witten differential} on it is a section $\bla_{SW}\in H^0(\cB,\cV^1)$, for which the $\cO_\cB$-module homomorphism
	
	\[
		T_\cB \longrightarrow \cV^0, \ v\longmapsto \nabla_v \bla_{SW}
	\]
	 factors through  an \emph{isomorphism}
	 \begin{equation}\label{isom_SW}
		T_\cB \simeq  \cV^1.
	\end{equation}

		Given such data, we obtain a refinement of the  weight-1 filtration
		\[
			\cV^1 \subseteq \cV^0
		\]
		to a  weight-3 filtration
		\[
			 \underbrace{\cF^3}_{=\bla_{SW}\cO_\cB} \subseteq  \underbrace{\ \cF^2 \ }_{\cV^1 }\subseteq \underbrace{\cF^1}_{=(\cF^3)^\perp}\subseteq 
			 \underbrace{\cF^0}_{\cV^0}
		\]

		 For links to  projective special K\"ahler geometry (``$N=2$ supergravity'') and weight-3 VHS, satisfying the Calabi--Yau condition, 
		 one can check \cite{hertling_hoev_posthum}[\S 4, \S 8.3].
		
		 Furthermore, given such data, there is an associated fibration of complex tori  $\cJ:= \cV/(\cV^1+\cV_\ZZ)\to \cB$,
		 whose vertical bundle is $\textrm{Vert}=\cV/\cV^1$. The polarisation $S$ gives rise to an isomorphism
		 $\textrm{Vert}\simeq (\cV^1)^\vee$, and hence $\bla_{SW}$ induces, by composition with the dual of its defining isomorphism $T_\cB\simeq \cV^1$,
		 an isomorphism $i_{\bla}: \textrm{Vert}\to T^\vee_\cB$. Such an isomorphism is also induced by a choice of symplectic form on $\cJ$.
		There is unique symplectic form $\omega_{\bla}$ on $\cJ$, which induces $i_{\bla}$ and such that the $0$-section is Lagrangian.

      We next recall the construction of the family of cameral covers over the Hitchin base, and introduce a weight-1 VHS with a Seiberg--Witten differential
      on it.

      First, we fix the following data:
	  \begin{itemize}
	  	\item A simple  complex Lie group $G$ of rank $l$,  together with a choice of Borel and Cartan subgroups $T\subset B\subset G$. We denote by
	  	$ \ft\subset \fb \subset \fg $ the respective Lie algebras and  by $W$ the corresponding Weyl group.

	  	\item A compact (connected) Riemann surface $X$ of genus $g\geq 2$ (or equivalently, a non-singular proper algebraic curve over $\CC$).
	  	We do not need to fix a particular projective embedding of $X$.
	  \end{itemize}
     Additionally, we choose:
	  \begin{itemize}
	  	\item Homogeneous generators $I_1,\ldots, I_l$ of the ring $\CC[\ft]^W\subset \CC[\ft]$. We write $d_k=\deg I_k$.
	  	\item Simple (positive) roots $\{\alpha_1,\ldots, \alpha_l\}$. 
	  \end{itemize}

	  These additional  choices are not necessary for the entire discussion, but are needed for the explicit calculation in
	  Theorem A.
	  
	  Two explicit examples of invariant polynomials -- for $SL_3(\CC)$ and $G_2$ -- are given in Equations (\ref{sl3}) and (\ref{G2_polys}), respectively.
	  
	  Notice that while $\ft/W$ is \emph{a priori} just a cone, the choice of generators $\{I_k\}$  allows us to identify it with $\CC^l$.
	  Notice also that we may   interpret $\{I_k\}$ as elements of $\CC[\fg]^G$, via Chevalley's theorem.
	  
	  The chosen simple roots determine an isomorphism $\ft\simeq \CC^l$, $v\mapsto (\alpha_1(v),\ldots, \alpha_l(v))$, 
	  using which we further identify
	  $\chi:\ft\to \ft/W$ with a finite map $\bI: \CC^l\to \CC^l$.
	  We may abuse the notation for these maps, e.g., write  $\chi=(I_1,\ldots, I_l)$ instead of $\bI$, etc.

	  We proceed by constructing from these data two rank-$l$ vector bundles on $X$.
	  The first one is $\ft\ctimes K_X \simeq K_X^{\oplus l}$, whose total space will be denoted by $M$: 
	  \[
	      M=\tot \ft\ctimes K_X.
	  \]
	  The group $W$ acts (fibrewise, via its action on $\ft$) on $M$. 
	  The resulting quotient $\bdU$ is \emph{a priori} just a cone bundle, but the choice of $\{I_k\}$ allows us to give it the
	  structure of a vector bundle of rank $l$:
	  \begin{equation}\label{U_isom}
	  	\bdU= \ft\ctimes K_X/W \simeq \bigoplus_{k=1}^l K_X^{d_k}.
	  \end{equation}
	  We can also think of $\bdU\backslash\{0\}$ as the $\CC^\times$-bundle with fibre $\ft/W$, 
	  associated to the $\CC^\times$-bundle $K_X\backslash\{0\}$.

	  The morphism $\chi:\ft\to \ft/W$ induces a morphism $\chi:M\to \tot \bdU$ of $X$-varieties (not of vector bundles!):
	  \begin{equation}\label{cover_MU}
	  	\xymatrix{M= \tot \ft\ctimes K_X \ar[rr]^-{\chi=(I_1,\ldots,I_l)}\ar[dr]^-{\pi}& & \tot \bdU\ar[dl]\\
	  	&X&\\}.
	  \end{equation}

	  We  write $\cB$ for the Hitchin base -- the space of global sections of $\bdU$:
	  \[
	  	\cB := H^0(X,\bdU)\simeq H^0\left(X, \bigoplus_{k=1}^l K_X^{d_k}\right) \simeq \CC^{\dim G(g-1)}.
	  \]

	  Any $b\in\cB$ determines a $W$-cover $p_b:\widetilde{X}_b\to X$ as the pullback of $\chi:M\to \tot \bdU$ via (the evaluation map of) the section $b$:
	  \[
	  	\xymatrix{ \widetilde{X}_b \ar[d]^-{p_b}\ar[r]& \tot \ft\ctimes K_X= M\ar[d]^-{\chi}\\
				X\ar[r]^-{ev_b}& \tot \ft\ctimes K_X/W=\tot \bdU\\}
	  \]
	  This $W$-cover is called \emph{the cameral cover of $X$ (corresponding to $b$)}. We may occasionally write $p:\widetilde{X}\to X$
	  if the point $b\in\cB$ is fixed or understood.
	  
	  By construction $\widetilde{X}_b$ is a closed subscheme of $M$ that can be singular or non-reduced.
	  The cameral cover $\widetilde{X}_b\subset M$ inherits from $M$ a $W$-action (and thus has lots of automorphisms). For a generic choice of $b$ it is a
	  non-singular ramified
	  Galois $W$-cover with simple ramification. We write  $\scB\subseteq \cB$ for the open set of generic cameral covers. 
	  
	   The vector bundle $p_b^\ast \bdU$ is in fact isomorphic to $N_{\widetilde{X}_b/M}$, the normal bundle of $\widetilde{X}_b\subseteq M$,
	  see Section \ref{emb_cam_curve}.
	  
	  \begin{Ex}
	  	Let $G=SL_2(\CC)$. Then $W=\ZZ/2\ZZ$, $\bdU\simeq K_X^2$,  	$\cB= H^0(X,K^2_X)\simeq \CC^{3g-3}$ and $p_b: \widetilde{X}_b\to X$ is a 
	  	$2:1$ cover. The open set $\scB\subseteq \cB$ consists of quadratic differentials with simple roots. For $b\in \scB$, the cover has  genus $g(\widetilde{X}_b)=4g-3$.
	  \end{Ex}
	  \begin{Ex}
	  	Let $G=G_2$. Then $W=D_6$ (dihedral group of order $12$) and $\bdU\simeq K_X^2\oplus K_X^6$.
	  	Consequently
	  	$\cB= H^0(X,K_X^2)\oplus H^0(X,K_X^6)\simeq \CC^{14(g-1)}$. The cameral covers
	  	$p_b:\widetilde{X}_b\to X$ are $12:1$ covers, with $g(\widetilde{X}_b)= 84(g-1)+1$.
	  \end{Ex}

	    There is a weight-1 $\ZZ$-variation of Hodge structures $\cV^1\subseteq \cV^0$
	  over $\scB\subseteq \cB$, 
	  whose fibres are respectively $\cV^1_b=H^0(\widetilde{X}_b,\ft\ctimes K_{\widetilde{X}_b})^W$ and $\cV_b^0= H^1(\widetilde{X}_b,\ft)^W$.
	    Intrinsically, it is defined as follows. Let $\Lambda\subseteq \ft$ be the cocharacter lattice and $p:\cX\to \scB$ the universal cameral curve.
	    Let also $p_\ast^W$ be the $W$-invariant pushforward functor. Then we set $\cV_\ZZ= R^1p_\ast^W(\Lambda)$ and
	    $\cV:=\cV_\ZZ\otimes_\ZZ\cO_\scB\simeq R^1p_\ast^W(\ft\ctimes \Omega^\bullet_{\cX/\scB})$.
	    The  bundle $\cV^1=R^0p_\ast^W\left(\ft\ctimes \Omega^1_{\cX/\scB}\right)$, and the Hodge filtration is 
	    induced by the naive filtration $\Omega_{\cX/\scB}^{\bullet \geq 1}[-1]\subseteq \Omega^\bullet_{\cX/\scB}$.
	    The Gauss--Manin connection can be identified with the $d_1$ differential of the spectral sequence, induced by the
	    Koszul--Leray filtration on $\Omega^\bullet_{\cX}$.
	   The polarisation pairing $S$ is given by
	   	$S_b(\alpha,\beta)=\left\langle\alpha\cup\beta, \left[\widetilde{X}_b\right]\right\rangle$.
	   	For more details, see section \ref{hhp} and the references therein, as well as \cite{hertling_hoev_posthum}[8.1]
	   	and \cite{don-pan}.

	  On $M$ there is a canonical $\ft$-valued Liouville form $\lambda$,
	    see section \ref{Liouville}.
	  The  Liouville form $\lambda$ determines  a Seiberg--Witten differential, 
	  $\bla_{SW}\in \Gamma(\scB,\cV^1)$, via
	   $\bla_{SW}(b)=\left. \lambda\right|_{\widetilde{X}_b}$, and, as in (\ref{isom_SW}), we have that
	    the map
	  \begin{equation}\label{key_isom_GM}
	  	\cB= T_{b}\scB\ni \bg \longmapsto \left(\nabla^{GM}_\bg \bla_{SW}\right)_b \in \cV^0_b
	  \end{equation}
	  factors through an isomorphism $T_{b}\scB\simeq \cV^1_b$, i.e.,

	  \begin{equation}\label{key_isom}
	  	\xymatrix@1{\left(\nabla^{GM}\bla_{SW}\right)_b: &  H^0\left(X, \bigoplus_{k=1}^l K_X^{d_k}\right)\ar[r]^-{\simeq} & H^0(\widetilde{X}_b,\ft\ctimes K_{\widetilde{X}_b})^W.}
	  \end{equation}

	  In  \cite{hurtubise_markman_rk2}[Proposition 2.11],
	  an   isomorphism with the same domain and codomain as in (\ref{key_isom}) is described as
	  the composition of
	  pullback on global sections (by $\pi$),  contraction with $\omega$ and restriction to $\widetilde{X}_b$,
	  see also Proposition \ref{isom_general_prop1}. 
	  In \cite{hertling_hoev_posthum}[Proposition 8.2] it is shown, using a hypercohomology calculation, that the isomorphism described
	  by Hurtubise and Markman coincides with the isomorphism  (\ref{key_isom_GM}).
	  Some of the above relations  for $G=SL_2$ are discussed in  \cite{douady_hubbard_strebel}[Proposition 1], see also 
	  \cite{msww_asymp}[Eq.(3)].

	  The above isomorphism can also be considered from an integrable systems viewpoint.  
	  Indeed, consider the  universal family of generic cameral curves $p:\cX\to \scB\subseteq \cB$.
	  The relative Prym fibration $\Prym_{\cX/\scB}\to \scB$ is in fact an
	  algebraic completely integrable
	  system. The fibre $\Prym_{\widetilde{X}_b}$ over  $b\in \scB$ is an abelian variety, whose tangent space is Serre dual to $H^0(\widetilde{X}_b,\ft\ctimes K_{\widetilde{X}_b})^W$,
	  the right hand side of (\ref{key_isom}). 
	  The isomorphism (\ref{key_isom}) actually amounts to  lifting a tangent vector in $T_{\cB,b}$ to a vector field
	  along the fibre $\Prym_{\widetilde{X}_b}$ and then 
	  pairing it with the symplectic form on the Prym fibration. This is the viewpoint, taken, e.g., by Hurtubise and Markman.

	  Our goal in this note is to provide an
	   explicit and global (on $X$ and $\widetilde{X}_b$) description of (\ref{key_isom}) in terms of Lie theory and the covering
	  $p_b:\widetilde{X}_b\to X$.

	  The simplest case, that of $G=SL_2$, is given in Example \ref{isom_sl2}, where we show that 
	  Equation (\ref{key_isom}) specialises to 

	  	\[
	  		H^0(X,K_X^2)\ni \bg\longmapsto \left(\nabla^{GM}_\bg \bla_{SW}\right)_b= \left. \frac{p^\ast \bg}{2\bal^2}\bla_{SW}\right|_{\widetilde{X}_b}\in 
	  		H^0(\widetilde{X}_b, K_{\widetilde{X}_b})^{\ZZ/2}_{-},
	  	\]
	  	where $\bal^2\in H^0(M,\pi^\ast K_X^2)$ is the tautological section and $\bla_{SW}$ is the Liouville (Seiberg--Witten) form.
	  	The expression on the right hand side can in fact also be rewritten as $-\frac{p^\ast g}{2\lambda}$, and in this form
	  	it coincides  (up to scaling
	  	factors) with  \cite{msww_asymp}[(4)], who  reference Douady--Hubbard \cite{douady_hubbard_strebel}[\S 2].

	  	Our main result is a general formula for $\nabla^{GM}_\bg\bla_{SW}$ for the case of an arbitrary (complex, simple) group $G$.
	  	
	  	Let  $D\bI$ be the Jacobi matrix of the adjoint quotient  $\bI=(I_1,\ldots, I_l):\CC^l\to\CC^l$ and  $\io$ the
		 natural algebra  homomorphism from
		$ \sym(\ft^\vee)$ into $ H^0\left(M, \bigoplus_{n\geq 0}\pi^\ast K_X^n\right)$, 
		introduced in Equation (\ref{making_bold}). 
		Finally, $\bal_i=\io(\alpha_i)$ and  $\lambda_i=e_i\otimes\bal_i$, where $\{e_i\}$ is the basis of $\ft$, dual to $\{\alpha_i\}$.
		In this notation, the Liouville form is $\bla_{SW}=\sum_i\lambda_i$.
	  
	  \begin{thma}\label{theorem_a}
		Once the
		  main and  additional data are chosen, the isomorphism (\ref{key_isom}) $\nabla^{GM}\bla_{SW}$
		  maps $\bg\in T_b\scB=H^0(X,\bigoplus_i K_X^{d_i})$ to  the section
		  \begin{equation}\label{main_formula}
		  	\left(\nabla^{GM}_\bg \bla_{SW}\right)_b = - \left. \sum_{i=1}^l \frac{\left( \io (D\bI)^{-1} \cdot \pi^{\ast}\bg\right)_i}{\bal_i}\lambda_i \right|_{\widetilde{X}_b}
		  		=
		  	 - \left. \io(D\bI)^{-1}\cdot \pi^\ast \bg\right|_{\widetilde{X}_b}.
		  \end{equation}
	\end{thma}
	
	  In particular, for  $l=2$ we have that
	  \[
	  	\left(\nabla^{GM}\bla_{SW}\right)_b:\ H^0(X, K_X^{d_1}\oplus K_X^{d_2}) \longrightarrow H^0(\widetilde{X}_b, K_{\widetilde{X}_b}^{\oplus 2})^W
	  \]
	  sends $\bg=	\begin{bmatrix}
	  	      		g_1\\ g_2\\
	  	      	\end{bmatrix}$ to

	  \begin{equation}
	  	 \left( \nabla^{GM}_\bg\bla_{SW}\right)_b= -
	  	  \left.
	  	\left(
	      \frac{
	      \begin{vmatrix}
                    \pi^*g_1& \io \partial_2 I_1 \\
                     \pi^*g_2& \io \partial_2 I_2 \\
		    \end{vmatrix}}
	      {\bal_1  \det  \io D\bI}\lambda_1 + 
	      \frac{
	      \begin{vmatrix}
                    \io \partial_1 I_1& \pi^*g_1 \\
                     \io \partial_1 I_2&  \pi^*g_2 \\
		    \end{vmatrix}}
	      {\bal_2  \det \io D\bI}\lambda_2
	      \right)\right|_{\widetilde{X}_b}
	      =- \frac{1}{\det\io D\bI} \left. 
	      \begin{bmatrix}
		    \begin{vmatrix}
                    \pi^*g_1& \io \partial_2 I_1 \\
                     \pi^*g_2& \io \partial_2 I_2 \\
		    \end{vmatrix} \\ 
		    \\
		    \begin{vmatrix}
                    \io \partial_1 I_1& \pi^*g_1 \\
                     \io \partial_1 I_2&  \pi^*g_2 \\
		    \end{vmatrix}
	      \end{bmatrix}
	      \right|_{\widetilde{X}_b}.
	  \end{equation}

	  	Knowledge of $\bla_{SW}$ and $\nabla^{GM}\bla_{SW}$ is essential for describing various geometric structures on $\scB$. 
	  	We mention only two examples as an illustration.

	  	First,  for the Hitchin integrable system, the Donagi--Markman cubic (\cite{donagi_markman_cubic}), which is essentially
	  	the infinitesimal period map for the family of Hitchin Pryms, is given by the Balduzzi--Pantev formula \cite{balduzzi}[Theorem 1].
	  	If we consider the cubic as a global section $c$ of $\sym^3 T^\vee_{\scB}=\sym^3\cB^\vee\otimes \cO_{\scB}$, then
	  	the Balduzzi--Pantev formula states that the value of $c$
	  	  at $b\in \scB$ is
	  	\begin{equation}
  c_b(\bg_1,\bg_2,\bg_3)=\frac{1}{2}\sum_{m\in \textrm{Ram }p_b}\textrm{Res}^2_m
  \left(p_b^\ast\left. \frac{\cL_{\bg_1} \mathfrak{D}}{\mathfrak{D}}\right|_{\{b\}\times X}\left(\nabla^{GM}_{\bg_2}\bla_{SW}\right)_b\cup \left(\nabla^{GM}_{\bg_3}\bla_{SW}\right)_b   \right).
  \end{equation}
	  	Here $\mathfrak{D}$ is the discriminant (see also section \ref{examples}) and $\cL$ denotes Lie derivative.
	  	In our previous work \cite{ugo_peter_cubic}[Theorem A] we have shown that the Balduzzi--Pantev formula holds along
	  	the (good) symplectic leaves of the generalised Hitchin system.

	  	The  second example which is worth mentioning is the special K\"ahler metric $g_{SK}$ on $\scB$.
	  	 It is known that for the case of $G=SL_2(\CC)$, the special K\"ahler metric  is given by
	  \begin{equation}
	  	g_{SK}(\bg,\bg)_b=2\int_{\widetilde{X}_b}|\nabla^{GM}_\bg\bla_{SW}|^2,
	  \end{equation}
	   see \cite{fredrickson_exp}[2.40], \cite{msww_asymp}[\S 2.3], \cite{dumas_neitzke}.

	  	   We shall discuss additional applications of {\bf Theorem A} to  various aspects of the geometry of
		  $\scB$ in a forthcoming work.

	    \subsection*{Acknowledgements} P.D.~ thanks Tony Pantev for helpful discussions related to the project,  
	    the Simons Collaboration on Homological Mirror Symmetry for support and the University of Pennsylvania for its hospitality.

	    \section{Preliminaries}

	      \subsection{The Embedding of the Cameral Curve}\label{emb_cam_curve}
	  We are now going to work at a fixed point $b\in \cB$ (generic), and hence will write mostly $p:\widetilde{X}\to X$ for the cameral cover. 
	  To understand (\ref{key_isom}) we need to understand $K_{\widetilde{X}}$ and for that we need to know more about the normal bundle $N$ of the closed
	  embedding $\widetilde{X}\subseteq M$. This is not difficult, since $\widetilde{X}$ is in fact the zero locus of a section of  a vector bundle on $M$.
	  
	  First, notice that the morphism $\chi:M\to \bdU$ (see (\ref{cover_MU})) induces a tautological section $\sigma\in H^0(M, \pi^\ast \bdU)$ in a standard way, via
	  \begin{equation}\label{taut_section}
	  	\xymatrix{ 
	  	M\ar@/^/[drr]^-{\chi}\ar@/_/[ddr]^-{\textrm{id}}\ar@{.>}[dr]^-{\sigma} &					&\\
										       &\pi^\ast \bdU\ar[d]\ar[r]		&\bdU\ar[d]\\
												  &M\ar[r]^-{\pi}	& X\\},
	  \end{equation}
	  which on closed points is simply $\sigma(m)= (m, \chi(m))\in M\times_X \tot \bdU=\tot\pi^\ast \bdU$.
	  
	  Next, the adjunction morphism $\bdU\to \pi_\ast\pi^\ast \bdU$ induces on global sections the pullback map 
	  $\cB=H^0(X,\bdU)\to H^0(M, \pi^\ast \bdU)$, which we write  as $b\mapsto \pi^\ast b$. 
	  
	  Thus the cameral curve $\widetilde{X_b}$ is the zero locus
	  \begin{equation}\label{cam_eqn}
	  	\widetilde{X}_b=\textrm{zeros}(s_b),\ s_b=\sigma - \pi^\ast (b)\in H^0(M, \pi^\ast \bdU),
	  \end{equation}
	  i.e., is cut out by the equation(s) 
	  \begin{equation}
	      \chi(m)=b(\pi(m))
	  \end{equation}
	  in $M=\tot \ft\ctimes K_X$. Having fixed  basic invariant polynomials $\{I_k\}$, 
	  and hence an isomorphism $\bdU\simeq \bigoplus_{k=1}^l K_X^{d_k}$, we can express this as the system of equations
	  \begin{equation}\label{cameral_global}
	  	\left|
		      \begin{array}{l}
		      	I_1(m)= b_1(\pi(m))\\
		      	\vdots \\
		      	I_l(m)= b_l(\pi(m))\\
		      \end{array}
	  	\right., 
	  \end{equation}
	  for $m\in M$, with  $b=(b_1,\ldots,b_l)\in \cB$ fixed. These are ``global'' equations and no choice of local trivialisation is used here: the $k$-th
	  equation takes values
	  in (the total space of) $K_X^{d_k}$. Another global description is given in Equation (\ref{cameral_global_2}).
	  
	  From Equation (\ref{cam_eqn})  follows 
	  
	  \begin{Prop}\label{claim_normal_bun}
	  The normal bundle of $\widetilde{X}_b\subseteq M$ is
	  \begin{equation}\label{normal_bundle} 
	  	N_{\widetilde{X}_b/M}\simeq p_b^\ast \bdU= \ft\ctimes p_b^\ast K_X/W \simeq \bigoplus_{k=1}^l p_b^\ast K_X^{d_k}.
	  \end{equation}
	  \end{Prop}
	  
	  \emph{Proof:}
	  While in general one uses the Koszul complex to compute the normal bundle, here we have
	  that both $\widetilde{X}_b$ and $M$ are smooth, and moreover, $\widetilde{X}_b$ is a complete intersection.
	  This case is handled by  a  standard geometric argument, given in, e.g.  \cite{eisenbud-harris_3264}[Proposition 6.15]. 
	  
	  The isomorphism $N_{\widetilde{X}_b/M}\simeq p_b^\ast \bdU$
	  is induced by the (vertical component of the) differential $ds_b: T_M\to s_b^\ast T_{\pi^\ast \bdU}$ of the section $s_b:M \to \tot\pi^\ast \bdU$. 

	  \qed

	  Similarly to the above argument, since $M$ is the total space of a vector bundle (namely $\ft\ctimes K_X$) on $X$, its tangent bundle $T_M$ is an extension of
	  $\pi^\ast T_X$ by $\ft\ctimes \pi^\ast K_X$.
	  Restricting to $\widetilde{X}$ and combining with the previous result,  one gets the diagram
	  \begin{equation}\label{display}
	  	\xymatrix{
				& 				& 0\ar[d]	 				& 			& \\
				& 				& T_{\widetilde{X}}\ar[d]			& 			& \\		
		     0\ar[r]	&\ft\ctimes p^\ast K_X\ar[r]	& \left. T_M\right|_{\widetilde{X}}\ar[r]\ar[d]	&p^\ast T_X\ar[r]	& 0\\
				&				&\ft\ctimes p^\ast K_X/W\ar[d]			&			& \\
				&				&0						&			&\\
			 }
	  \end{equation}

	Now,
	consider $\bg=(g_1,\ldots,g_l)\in T_{\cB,b}=\cB$, with $g_i\in H^0(X,K_X^{d_i})$. It determines
	  a $1$-parameter family of deformations of $\widetilde{X}_b$,  given
	  by the equation
	  \begin{equation}\label{defo_family}
	  	\chi(m)= b(\pi(m))+\epsilon \bg(\pi(m)),
	  \end{equation}
	  that is, $\left\{ \widetilde{X}_{b+\epsilon g}\right\}_\epsilon$.
	  For $\epsilon$ in a sufficiently small disk  $\Delta_\rho\subseteq \CC$ the section  $b+\epsilon \bg\in \cB$ remains generic -- which we assume to be the case from
	  now on. The total space of the $1$-parameter family is cut out  in $M\times \Delta_\rho$ by the Equation (\ref{defo_family}).

	  The section $\bg$ determines a section of $N_{\widetilde{X_b}/M}= \ft\ctimes p^\ast_b K_X/W$, namely, $p_b^\ast \bg$.

	      \subsection{Local Description}\label{local_1}
	
	It is not hard to describe the objects from the previous section  in local coordinates.
	A choice of  a local (analytic) chart $\psi$ on $X$, identifying an open $U\subseteq X$ with a disk $\Delta\subseteq \CC$,
	 determines a local trivialisation of $K_X$
	and a compatible
	bundle chart $\phi$ on $M$, identifying $M_U= \pi^{-1}(U)\to U$ with $\textrm{pr}_1:\Delta\times \ft\to\Delta$,
	as usual:
	\begin{equation}\label{loc_coord}
		\xymatrix{M_U=\pi^{-1}(U)\ar[r]^-{\phi}\ar[d]_-{\pi} & \Delta\times \ft\ar[d]^-{\textrm{pr}_1}\\
			X\supseteq  U\ar[r]^-{\psi}& \Delta\subseteq \CC}.		  
	\end{equation}
	Such a  local chart   determines a trivialisation of
	$K_X^{d_i}$ over $U$ and hence a section $b_i\in H^0(X,K_X^{d_i})$ is represented locally as $(\psi^{-1})^\ast b_i= \beta_i(z) dz^{\otimes d_i}$ on $U$, 
	where $\beta_i:\Delta\to \CC$
	is a holomorphic function.

	Using the   simple roots as a basis  for $\ft\simeq \CC^l$, we identify  $\widetilde{X}_U=p_b^{-1}(U)$  (via $\phi$) with
	the set of solutions of $\bI(\alpha_1,\ldots,\alpha_l)=\bbe(z)$
	  for $(z,\underline{\alpha})\in \Delta\times \CC^l$, giving a local version of Equation (\ref{cameral_global}).

	  Next, the trivialisations of $K_X^{d_i}$ ($i=1\ldots l$) and the choice of roots provide  an induced trivialisation 
	  $\left. T_{M_U}\right|_{\widetilde{X}_U}$ and
	  \begin{equation}
	  	\left. (\psi^{-1})^*T_{M_U}\right|_{\widetilde{X}_U}= \left(\CC\oplus \ft\right)\ctimes \cO_{\phi(\widetilde{X}_U)}\simeq  \cO_{\phi(\widetilde{X}_U)}\left\langle \frac{\partial}{\partial z}, \frac{\partial}{\partial \alpha_1},
	  	\ldots , \frac{\partial}{\partial \alpha_l}\right\rangle
	  \end{equation}
	  and, consequently, a local description of the diagram (\ref{display}):
	  \begin{equation}\label{display_loc}
	  	\xymatrix{
				& 				  & 0\ar[d]	 				& 			& \\
				& 				  & (\psi^{-1})^\ast T_{\widetilde{X}_U}\ar[d]			& 			& \\		
		     0\ar[r]	&\ft\ctimes \cO_{\phi(\widetilde{X}_U)} \ar[r]^-{ \begin{pmatrix}0\\ 1\\ \end{pmatrix} }	&  \left(\CC\oplus \ft\right)\ctimes \cO_{\phi(\widetilde{X}_U)} \ar[r]^-{ \begin{pmatrix} 1& 0 \end{pmatrix} }\ar[d]^-{ \begin{pmatrix} -\bbe'& D\bI \end{pmatrix} }	& \cO_{\phi(\widetilde{X}_U)} \ar[r]	& 0\\
				&				  &\CC^l\ctimes \cO_{\phi(\widetilde{X}_U)}\ar[d]			&			& \\
				&				  &0						&			&\\
			 }.
	  \end{equation}

	  Here the bottom vertical map is, in more detail, 
	  \begin{equation}
	  	 \begin{pmatrix} -\bbe'& D\bI \end{pmatrix}= 
	  	 \left( 
		      \begin{array}{rccc}
		      	-\beta_1'& \partial_1 I_1& \ldots & \partial_l I_1\\
		      	\vdots   & \vdots	 & \vdots & \vdots \\
		      	-\beta_l'& \partial_1 I_l& \ldots & \partial_l I_l\\
		      \end{array} 
	  	 \right)\in \textrm{Mat}_{l\times (l+1)}\left(\Gamma\left(\cO_{\phi(\widetilde{X}_U)}\right)\right),
	  \end{equation}
	  having  rank $l$ everywhere on $\widetilde{X}_U$, under the assumption that
	   $b=(b_1,\ldots,b_l)\in \cB$ is generic.
	  This is  the matrix of the map $\textrm{pr}_2\circ ds$ from Proposition \ref{claim_normal_bun}.
	  We write    $D\bI$ or $D\chi$ for the Jacobi matrix of $\bI=(I_1,\ldots, I_l):\CC^l\to\CC^l$.

	  Finally, given a tangent vector $\bg=(g_1,\ldots,g_l)\in T_{\cB,b}=\cB$, with $(\psi^{-1})^* g_i= \gamma_i(z)dz^{\otimes d_i}$ on $U$,
	  the corresponding
	   $1$-parameter (analytic) family of deformations of $\widetilde{X}_b$ is  cut out locally (in $\Delta \times \CC^l\times \Delta_\rho$) by
	  $\bI(\underline{\alpha})=\bbe(z) + \epsilon \bga (z)$,
	  where $\Delta_\rho\subseteq \CC$ is as before.

	 We may occasionally suppress the pullbacks by $\phi$ and $\psi$, except for the cases when there is a risk of confusion, as when
	discussing (co)roots and some associated objects.

	\subsection{Objects, Associated with Roots}\label{symp_form}

	Any linear map $\alpha\in\ft^\vee=\mhom(\ft,\CC)$ determines, by extension of scalars,  a vector bundle homomorphism $\ft\ctimes K_X\to K_X$,
	denoted by the same letter.
	Hence, just as $\chi$ in Eq (\ref{taut_section}), such an
	$\alpha$ determines a tautological section $\bal\in H^0(M, \pi^\ast K_X)$, which on (closed) points  maps $m\in M$ to 
	   $ \bal(m)=(m, \alpha(m))\in M\times_X \tot K_X$.
	Furthermore, 
	 restricting  $\bal$  to $\widetilde{X}\subset M$ gives  a section
	$\bal_{\widetilde{X}}\in H^0(\widetilde{X},p^\ast K_X)$. Occasionally, we  suppress the subscript
	$\widetilde{X}$, i.e., the restriction.

		The section $\bal$  vanishes along a  ``hyperplane divisor'' 
		      $\tot \left(\ker \alpha\ctimes K_X\right) \subseteq M$,
		  a rank-$(l-1)$ subbundle of $\ft\ctimes K_X$.
	    	The respective restrictions
	$\bal_{i \widetilde{X}}$ (of sections arising from roots) vanish along divisors $D_{\alpha_i}$ in $\widetilde{X}$, which are
	the ramification divisors of $p:\widetilde{X}\to X$.

	If we choose a local chart $(U,\psi)$ and $\phi:M_U\simeq \Delta\times \ft$,  as in  (\ref{loc_coord}),
	$\bal$
	is represented by $(z,u)\mapsto \alpha(u)dz$, 
	where $\alpha(u)=\langle \alpha,u \rangle$ is the natural pairing between $\ft$ and $\ft^\vee$.
	If we further identify the preimage of $\pi^{-1}(U)$ in $\tot\pi^\ast K_X\to M$ with $\Delta\times\ft\times \CC$, 
	via $\phi$ and a trivialisation of $K_X$,
	then
	the evaluation map
	of $\bal$ is represented  by
	\[
		\Delta\times\ft\ni	(z,u) \longmapsto (z,u,\alpha(u)) \in \Delta\times\ft\times \CC.
	\]

	The linear functional  $\alpha\in \ft^\vee$ determines a function on $\Delta\times \ft$,  that we may 
	denote  $\textrm{pr}_2^\ast\alpha$ if the distinction from $\alpha$ is important. 
	Furthermore, given the choice of $\phi$, we may consider $\alpha$ (or rather, $\textrm{pr}_2^\ast\alpha$)
	a  function 
	$\phi^\ast \alpha\in \cO_{M_U}(M_U)$
	on $ M_U$.
	Consequently, upon restriction
	to $\widetilde{X}_U$, we get a local function $\phi^\ast \alpha\in \cO_{\widetilde{X}_U}(\widetilde{X}_U)$
	on the cameral curve.
	Of course, one should really write  $\left. \phi^\ast\textrm{pr}_2^\ast\alpha\right|_{\widetilde{X}_U}$ here.

	 The distinction between the various objects associated to a root $\alpha_i$ becomes important when one considers their
		differentials. Since  $\pi^\ast K_X\subseteq \Omega^1_M$,
			$d\bal_i\in \Omega^2_M(M)$. At the same time, 
		$d\phi^\ast \alpha_i\in \Omega^1_{M_U}(M_U)$ and
			$d(\textrm{pr}_2^\ast\alpha_i)\in \Omega^1(\Delta\times \ft)$.
			 Naturally, we are going to write $d\alpha_i$ for the penultimate expression, so the distinction between
			 $d\alpha_i$ and $d\bal_i$ is essential.
		Finally, we keep in mind that  $d\alpha_i=\alpha_i\in \mhom(\ft,\CC)$, as with any linear map.

		 The assignment $\alpha_i\mapsto \bal_i$ determines an (injective) $\CC$-algebra homomorphism
		\begin{equation}
			\io: \sym(\bt^\vee)\hookr H^0\left(M, \bigoplus_{n\geq 0}\pi^\ast K_X^n\right)
		\end{equation}
		and, consequently, a homomorphism
		\begin{equation}\label{making_bold}
			 \mend(\CC^l)\otimes \sym(\ft^\vee)\hookr \mend(\CC^l) \otimes H^0\left(M, \bigoplus_{n\geq 0}\pi^\ast K_X^n\right),
		\end{equation}
		denoted by $\io$ as well.
		Given a $\sym(\ft^\vee)$-valued endomorphism $A$ with non-zero determinant $\det A\in \sym(\ft^\vee)$, we write 
		$\io(A)^{-1}$ for the inverse of $\io(A)$ in the ring of $l\times l$ matrices with coefficients in the field of fractions
		$\textrm{Frac }H^0\left(M, \bigoplus_{n\geq 0}\pi^\ast K_X^n\right)$,
		and in fact, in 
		$\mend(\CC^l)\otimes H^0\left(M, \bigoplus_{n\geq 0}\pi^\ast K_X^n\right)\left[\frac{1}{\det \io(A)}\right]$.

		We can, more generally, rewrite the global equations for $\widetilde{X}_b$ as
	
	\begin{equation}\label{cameral_global_2}
	  	\left|
		      \begin{array}{l}
		      	I_1(\bal_1,\ldots,\bal_l)= \pi^\ast b_1 \\
		      	\vdots \\
		      	I_l(\bal_1,\ldots,\bal_l)= \pi^\ast b_l\\
		      \end{array}
	  	\right.
	  \end{equation}
	  that is, the linear system $\io(I_k)=\pi^\ast b_k$, $k=1\ldots l$.

	      \subsection{Liouville Form}\label{Liouville}
	      
	On $M=\tot \ft\ctimes K_X$ there is a $\ft$-valued 2-form $\omega\in H^0(M, \ft\ctimes \Omega^2_M)$. Probably the simplest way to introduce it is by setting
	\[
		\omega = - d\lambda,
	\]
	where $\lambda\in H^0(M,\ft\ctimes \pi^\ast K_X)\subseteq H^0(M,\ft\ctimes \Omega^1_M)$ is a tautological section, the ``$\ft$-valued Liouville form''.

	We  recall some explicit expressions for $\lambda$  -- although, as
	usual in symplectic geometry, there are various sign ambiguities in the possible definitions.      
	
	The chosen simple roots $\{\alpha_1,\ldots,\alpha_l\}$ form a basis of $\ft^\vee$, and we let $\{e_1,\ldots, e_l\}$ stand for the corresponding dual basis of $\ft$
	(consisting of fundamental coweights).

	One can then set  $\lambda_i= e_i\ctimes \bal_i$, a global section of $\ft\ctimes \pi^\ast K_X\subseteq \ft\ctimes \pi^\ast \Omega^1_M$, and
	write the Liouville form and the 2-form as
	\begin{equation}\label{liouville}
		\lambda= \sum_{i=1}^l \lambda_i=\begin{bmatrix}
		         	\bal_1\\ \bal_2\\ \vdots \\ \bal_l\\
		         \end{bmatrix},\ \omega = - \sum_{i=1}^l e_i\ctimes d\bal_{i}= \begin{bmatrix}
		        	-d\bal_1\\ -d\bal_2\\ \vdots\\ - d\bal_l\\
		        \end{bmatrix}.
	\end{equation}

	Finally, if we choose local coordinates as in eq  (\ref{loc_coord}),
	we obtain for the pullback of $\lambda$ and $\omega$  to $\Delta\times \ft$
	\[
		(\phi^{-1})^\ast   \lambda = \sum_{i=1}^n e_i\ctimes \alpha_i dz = 
				\begin{bmatrix}
	                           	\alpha_1 dz \\
	                           	\vdots\\
	                           	\alpha_l dz\\
				\end{bmatrix},\
				 (\phi^{-1})^\ast \omega = \sum_{i=1}^l e_i\ctimes dz\wedge d\alpha_i = 
	       \begin{bmatrix}
	                           	dz\wedge d\alpha_1\\
	                           	\vdots\\
	                           	dz\wedge d\alpha_l\\
				\end{bmatrix}.
	\]

	\section{Background: Two Results}
	\subsection{A Result of Hurtubise and Markman }\label{isom_general}
	We begin with the
	special case of a  result of Hurtubise and Markman \cite{hurtubise_markman_rk2}[Proposition 2.11] mentioned in the introduction.
	We spell out some of the details of their argument for this special case.

	\begin{Prop}\label{isom_general_prop1}
		For each generic $b\in \cB$, the pullback of global sections via $p_b$, followed by the isomorphism (\ref{normal_bundle}) and  contraction with  $\omega$  induces an isomorphism
		
		\begin{equation}\label{key_isom_2}
			\xymatrix@1{\beta: \cB= H^0(X, \bdU)\ar[r]^-{\simeq} & H^0(\widetilde{X}_b, p_b^\ast \bdU)^W \ar[r]^-{\simeq} &  H^0(\widetilde{X}_b, \ft\ctimes K_{\widetilde{X}_b})^W,  }
		\end{equation}

		or, using the choice of invariant polynomials $\{I_k\}$, an isomorphism

		\[
			 H^0\left(X, \bigoplus_{k=1}^l K_X^{d_k}\right)\simeq H^0\left(\widetilde{X}_b, \bigoplus_{k=1}^l p_b^\ast K_X^{d_k}\right)^W 
			 \simeq H^0\left(\widetilde{X}_b, K_{\widetilde{X}_b}^{\oplus l}\right)^W.
		\]

	\end{Prop}

	Thus, the isomorphism $\beta$ (\ref{key_isom_2}) is a composition of two maps.  The first one is pullback (adjunction) $\bg\longmapsto p_b^\ast \bg$, for $\bg\in \cB=H^0(X,\bdU)$. The second one 
	is the map on global sections, induced by the map of bundles
	
	\begin{equation}\label{bundle_map}
		N_{\widetilde{X}} \longrightarrow \ft\ctimes K_{\widetilde{X}}
	\end{equation}
	\[
		s\longmapsto \left. \omega(\widetilde{s},\ ) \right|_{\widetilde{X}},
	\]
	where $\widetilde{s}$ is a lift of $s$ to  a section of  $T_M$. One  may  denote this map simply by $\lrcorner\ \omega$ (contraction with $\omega$), but 
	should keep in mind the
	restriction to $\widetilde{X}$. 
	
	The proof of  Proposition \ref{isom_general_prop1} relies on a dimension count, combined with good understanding of
	the bundle map (\ref{bundle_map}) and  the induced map on fibres at $m\in \widetilde{X}$. For that,  the cases
	when $m$ is not a ramification point and when it is one should be considered separately.
	Notice that if $m$ is not a ramification point, then $T_{\widetilde{X},m}\nsubseteq \pi^{-1}(\pi(m))=\ft\ctimes K_{X,p(m)}$, while
	$T_{\widetilde{X},m}\subseteq \pi^{-1}(\pi(m))$
	if $m$ is a ramification point.
	
	So let us choose a point $m\in M$ and consider the fibre of $\pi:M\to X$, passing through $m$. We set
	  $L:= \pi^{-1}(\pi(m))=\ft\ctimes K_{X,\pi(m)}\subseteq M$, and write $N_L$ for the normal bundle of
	  the vector space $L\subseteq M$.

	Using the local description of $\omega$, we obtain that $\lrcorner\ \omega$ fits in the following diagram:
	
	\[
		\xymatrix{T_{L,m}\otimes N_{L,m}\ar@{^{(}->}[r] & \left. T_{M,m}\right|_{L}\otimes N_{L,m}\ar[r]^{\lrcorner \omega}\ar@{.>}[dr]&\ft\otimes \left. T^\vee_{M,m}\right|_L \otimes N_{L,m}\\
								&									       & \ft\otimes N_{L,m}^\vee\otimes  N_{L,m}= \ft\ar@{^{(}->}[u]\\
		}
	\]

	Since $N_L= T_{X,\pi(m)}\ctimes \cO_L$ and $T_L=L\ctimes\cO_L= \ft\ctimes K_{X,\pi(m)}\ctimes\cO_L$, there is a canonical
	trivialisation $T_L\otimes N_L=\ft\ctimes \cO_L$.
	
	Using  the normal sequence for $L\subseteq M$, one obtains:

	\begin{Lemma}\label{lemma0}
		 The map $\lrcorner\ \omega$ induces a trivialisation $T_L\otimes N_L\simeq_\omega \ft\ctimes\cO_L$, which coincides
		 up to sign with the canonical trivialisation. That is, 
		      \[
		\xymatrix{   T_L\otimes N_L \ar@{=}[d]_-{\textrm{can}}\ar[r]^-{\simeq_\omega} & \ft\ctimes \cO_L\\
		 \ft\ctimes \cO_L\ar[ur]^-{-\textrm{id}}& \\
		}.
		      \]
	\end{Lemma}

	Thus,  in particular,
	 $\lrcorner \omega$ induces, for any $m\in L$,   a $W$-equivariant isomorphism $T_{L,m}\otimes N_{L,m}\simeq \ft$.
	A similar result is stated,  in a much more general setup, in \cite{hurtubise_markman_rk2}[Theorem 2.8 (5)].

	\begin{Lemma}\label{lemma1}
		Consider a point  $m\in \widetilde{X}$ that is \emph{not} a ramification point of $p: \widetilde{X}\to X$. 
		The map on fibres, induced by the bundle map (\ref{bundle_map}) is   an isomorphism
		\[
			\xymatrix@1{\lrcorner\ \omega: & N_{\widetilde{X},m} \ar[r]^-{\simeq}&\  \ft\ctimes K_{\widetilde{X},m }.}
		\]
	\end{Lemma}

	 This is again a local calculation, using the explicit form of $\omega$. Notice that since $m$ is not a ramification point, the composition
	\[
		\xymatrix{  L=T_{L,m}\ar@{^{(}->}[r]  & T_{M,m}\ar@{->>}[r] & N_{\widetilde{X},m}\\}
	\]
	is an isomorphism. However, at ramification points
	the behaviour of $\lrcorner \omega$ is different.
	In fact, at such points  the map (\ref{bundle_map}) is \emph{not} an isomorphism of bundles if $l>1$, as is clear from the next Lemma.  
	
	\begin{Lemma}\label{lemma2}
		Let $m\in \widetilde{X}$ be a ramification point of $p:\widetilde{X}\to X$. Then
		\[
			\xymatrix{	N_{\widetilde{X},m}\ar[r]^-{\lrcorner\omega}\ar[dr] & \ft\ctimes K_{\widetilde{X},m}\simeq_\omega T_{L,m}\otimes N_{L,m}\otimes K_{\widetilde{X},m}\\
											    & T_{\widetilde{X},m}\otimes N_{L,m}\otimes K_{\widetilde{X},m}\ar@{^{(}->}[u]\\
				  }
		\]
		commutes.

	\end{Lemma}
	
	This result is shown by a local calculation, which in turn boils down to a linear-algebraic
	result, using the explicit form of $\omega$. It  is also stated in \cite{hurtubise_markman_rk2}[Lemma 2.10].

	\emph{Proof of Proposition \ref{isom_general_prop1}:}
	
	The map $\beta$, i.e., (\ref{key_isom_2}) is a composition of two maps, both of which are injective. Indeed, 
	$H^0(X,\bdU)\hookr H^0(\widetilde{X}, p^\ast\bdU)$, and the image is contained in $H^0(\widetilde{X}, p^\ast\bdU)^W$. 
	Furthermore, Lemma \ref{lemma1} implies that the map (\ref{bundle_map}) induces an injection on global sections,
	$H^0(\widetilde{X}, \ft\ctimes p^\ast K_X/W)\hookr H^0(\widetilde{X}, \ft\ctimes K_{\widetilde{X}})$, and it
	preserves $W$-invariant sections. Finally, by Serre duality (and the fact that taking duals commutes with taking invariants), we get
	\[
		\cB= H^0(X,\bdU)\hookr H^0(\widetilde{X}, p^\ast \bdU)^W\hookr H^0(\widetilde{X}, \ft\ctimes K_{\widetilde{X}})^W\simeq H^1\left(\widetilde{X},\ft\ctimes \cO_{\widetilde{X}}\right)^{W\vee}.
	\]
	But $ H^1\left(\widetilde{X},\ft\ctimes \cO_{\widetilde{X}}\right)^{W}$ is the tangent space of the generalised Prym variety, 
	and by the complete integrability of the
	Hitchin system, its dimension equals the dimension of the base $\cB$. 
	Hence  both  injections are isomorphisms, and so is their composition. 
	
	\qed

	\subsection{A Result of Hertling, Hoevenaars and Posthuma}\label{hhp}
	
	We introduced earlier a certain weight-1 $\ZZ$-VHS $(\cV, \nabla^{GM}, \cV_\ZZ, \cV^\bullet, S)$ on $\scB\subseteq \cB$.

	The bundle of lattices $\cV_\ZZ$ was defined as $R^1p_*^W(\Lambda)$, where $p: \cX\to \scB$ is the universal
	cameral cover, and the vector bundle $\cV=\cV_\ZZ\otimes_\ZZ\cO_{\scB}\simeq R^1p_*^W(\ft\ctimes p^*\cO_{\scB})$.
	The relative holomorphic Poincar\'e Lemma  gives a 
	quasi-isomorphism $p^{-1}\cO_{\scB}\simeq_{quis}\Omega^\bullet_{\cX/\scB}$,
	leading to $\cV\simeq R^1p_*^W(\ft\ctimes \Omega^\bullet_{\cX/\scB})$.
	
	The Hodge bundles, as for geometric VHS, are determined by the naive filtration of $\Omega^\bullet_{\cX/\scB}$,
	see \cite{voisin1}[\S 10.2] and \cite{voisin2}[\S 5.1]. In our case of weight one, 
	$\Omega^{\bullet \geq 1}_{\cX/\scB}[-1]\subseteq \Omega^\bullet_{\cX/\scB}$ determines a subbundle
	$\cV^1\subseteq \cV$, as $R^1p_*^W(\ft\ctimes \Omega^{\bullet\geq 1}_{\cX/\scB}[-1])\simeq R^0p_*^W(\ft\ctimes \Omega^1_{\cX/\scB})$.
	
	The Gauss--Manin connection $\nabla^{GM}:\cV\to \cV\otimes \Omega^1_{\scB}$  can be defined in either topological
	or holomorphic terms.
	 The topological description relies on Ehresmann's theorem, i.e., on  the $C^\infty$-local triviality of 
	 $p: \cX\to \scB$. In this case, the  homotopy-invariance of de Rham cohomology implies that
	$\cV_\ZZ$ is a locally constant sheaf  and $\nabla^{GM}$ can be described by a Cartan--Lie formula.
	For geometric VHS  this is described, e.g., in \cite{voisin1}[\S 8.2].
	
	The holomorphic description of $\nabla^{GM}$ is discussed in \cite{hertling_hoev_posthum}[\S 8], following
	\cite{katz-oda}, see also \cite{deligne_travaux_griffiths} and \cite{voisin2}[\S 5.1].
	 The Koszul--Leray
	filtration on $\Omega^\bullet_{\cX}$ gives rise to a spectral sequence, for which
	$(E_1^{\bullet,0},d_1)$ is identified with $(\Omega^\bullet_{\scB}(\cV),\nabla^{GM})$.
	
	One has the following result.
	
	\begin{thm}[\cite{hertling_hoev_posthum}, Proposition 8.2]\label{hhp_thm}
	      The isomorphisms $\nabla^{GM}\bla_{SW}$ and $\beta$ (\ref{key_isom_2}) coincide. That is, 
	      \[
	      	\nabla^{GM}_\xi\bla_{SW}=\beta(\xi),\ 
	      \]
	    for all tangent vectors $\xi\in T_b\scB$ and all $b\in\scB$.
		
	\end{thm}

	      The result is proved by an explicit hypercohomology calculation, using  the \v{C}ech resolution
	      of the relative de Rham complex $(\Omega^\bullet_{\cX/\scB},d)$.

	\section{Proof of  Theorem A}
	
	We now turn to the proof of our main result, {\bf Theorem A}.
	Recall that in the statement of the theorem  we use the algebra homomorphism $\io$ from Equation (\ref{making_bold}),   
	so $\io (D\bI)^{-1}$ is a global meromorphic section of $\send \left( \CC^l \ctimes \bigoplus_{k\geq 0}\pi^\ast K_X^k \right)$
	with poles  along the zeros of $\det\io D\bI$.
	That is, 
	the homogeneous polynomials $\partial_i I_j\in \sym^{d_j-1}(\ft^\vee)$ are considered as global sections of  $\pi_b^\ast K_X^{d_j-1}$, or, after
	restriction to $\widetilde{X}_b$, as sections of $p_b^\ast K_X^{d_j-1}$.

	 Using Cramer's formula and the fact that $\io$ is an algebra homomorphism, we can rewrite the right side of  (\ref{main_formula}) as
	 a linear combination of (restrictions of) $\lambda_i$ with coefficients of the kind
	\[
		\left.   \frac{\det\left[\io\partial_1 \bI,\ldots, p_b^\ast \bg, \ldots,\io \partial_l \bI  \right]}{\bal_i \det \io D\bI} \right|_{\widetilde{X}_b}\in \cK(\widetilde{X}_b),
	\]
	i.e.,  global meromorphic functions on $\widetilde{X}_b$, since both the numerator and the denominator belong to
	$H^0\left(\widetilde{X}_b, p_b^\ast K_X^{\sum_i d_i - l+ 1} \right)$.

	Now we prove Theorem A.
	Let us fix $\bg\in T_b\cB=H^0(X,\bigoplus_i K_X^{d_i})$ and denote by $\bs\in H^0(\widetilde{X}_b, \ft\ctimes K_{\widetilde{X}_b})$ the image of $\bg$ under 
	the isomorphism (\ref{key_isom}).
	Let us also denote by $\widetilde{\bs}$ the section from the right hand side of Equation (\ref{main_formula}), i.e., 
	\[
		\widetilde{\bs} = - \left. \sum_{i=1}^l \frac{\left( \io (D\bI)^{-1} \cdot \pi^{\ast}\bg\right)_i}{\bal_i}\lambda_i \right|_{\widetilde{X}_b}.
	\]
	This is a meromorphic section of $\ft\ctimes K_{\widetilde{X}_b}$  with poles at most along the ramification of $p_b:\widetilde{X}_b\to X$.
	We are going to prove that $\bs=\widetilde{\bs}$. We use Theorem \ref{hhp_thm} and the representation of the 
	isomorphism $\beta$ from Equation (\ref{key_isom_2}) is a composition of two maps.
	
	As a first step, we show that 
	\begin{equation}\label{Step_1}
		\left. \bs\right|_{\widetilde{X}_b\backslash \textrm{Ram}(p_b) } = \left. \widetilde{\bs}\right|_{\widetilde{X}_b\backslash \textrm{Ram}(p_b) }
	\end{equation}

	For that we restrict the cameral cover to the complements of the ramification  and branch divisors
	\[p_b: \widetilde{X}_b\backslash \textrm{Ram}(p_b) \longrightarrow X\backslash \textrm{Bra}(p_b)\]
	and choose  $U\subseteq X\backslash \textrm{Bra}(p_b)$, 
	  biholomorphic to an open disk (via $\psi:U\to \Delta$).
	  In this case, $\widetilde{X}_U \subseteq \widetilde{X}\cap \left(\det \io D\bI \neq 0\right)$  has $|W|$ (analytic) connected components,
	  each isomorphic to $U$, labelled by the different
	   Weyl chambers
	  \[
	  	\widetilde{X}_U=\widetilde{X}_U^1\coprod \ldots \coprod \widetilde{X}_U^{|W|}.
	  \]
	We choose (an  analytic)  local coordinate $z$ on $U$
	and use $z$ (i.e., its pullback $p_b^\ast z$) as a coordinate on $\widetilde{X}_U\subseteq \widetilde{X}_b\backslash \textrm{Ram}(p_b)$.

	  Then, setting $\bga$ for the coordinate vector of $p_b^\ast \bg$,
	  \[
	  	(\phi^{-1})^\ast\bga= 
			\begin{pmatrix}
	  	      	-\bbe'& \left. D\bI\right|_{\phi(\widetilde{X}_U)}\\
	  	      \end{pmatrix}
		      \begin{pmatrix}
		      	0\\
		      	\left.(D\bI)\right|^{-1}_{\phi(\widetilde{X}_U)}((\phi^{-1})^\ast\bga)\\
		      \end{pmatrix}\in \cO_\Delta^{\oplus l}(\Delta),
	  \]
	  i.e., we obtain a lift $\widetilde{\bga}$ of $\bga$ 
	  \begin{equation}\label{lift_1}
	  	 (\phi^{-1})^\ast \widetilde{\bga}= 
	  	 \sum_{i=1}^l (\left. (D\bI)^{-1}\right|_{\phi(\widetilde{X}_U)}(\phi^{-1})^\ast\bga)_i \frac{\partial}{\partial \alpha_i}\in \Gamma(\phi(\widetilde{X}_U),\left. (\phi^{-1})^*T_M\right|_{\widetilde{X}_U}).
	  \end{equation}

	  Note that the  expression for $\widetilde{\bga}$ is  well-defined on $\widetilde{X}_U$: away from ramification,
	  we can solve locally-analytically  for
	  $\alpha_i$ in terms of $z$, so $(D\bI)^{-1}$, when restricted to a connected component of $\phi(\widetilde{X}_U)$, is actually a
	  section of $\mend(\CC^l)\ctimes \cO^{an}_{\Delta}$.

	Then, using the lift $\widetilde{\bg}$ from Equation (\ref{lift_1}), we obtain
	\begin{multline*}
		\left. \widetilde{\bg}\lrcorner \omega\right|_{\widetilde{X}_U} = - \phi^\ast \left( \sum_{i=1}^l  \left(\left.(D\bI)\right|_{\phi(\widetilde{X}_U)}^{-1}\cdot \bga\right)_i \otimes e_i\otimes [dz] \right) \\ 
		- \sum_{i=1}^l  \left. \frac{\phi^\ast\left(\left.(D\bI)\right|_{\phi(\widetilde{X}_U)}^{-1}\cdot \bga\right)_i}{\alpha_i} \lambda_i\right|_{\widetilde{X}_U}, 
	\end{multline*}
	as $(\phi^{-1})^\ast\lambda_i= \alpha_i [dz]\otimes e_i$. 
	We write $[dz]$ rather than $dz$ since the cotangent sheaf of $\widetilde{X}_U$ is a quotient of  $\left. \Omega^1_{M_U}\right|_{\widetilde{X}_U}$.
	This is precisely the expression for $\widetilde{\bs}$ from Equation (\ref{main_formula}), written locally.
	
	Having shown (\ref{Step_1}),    we now  note that the sheaf of meromorphic sections of a   holomorphic vector  bundle on a smooth curve is trivial (see e.g.~\cite{GunningVB}, p.76; see also \cite{Stacks}, Lemma 31.25.3). As two meromorphic functions that coincide away from a finite set of points are equal, equation  (\ref{Step_1}) shows that $\bs = \widetilde{\bs}$. 
%
%
	Since the two sections $\bs$ and $\widetilde{\bs}$ are equal, and $\bs$ is known to be $W$-invariant, so is $\widetilde{\bs}$.
	\qed

	\section{Examples}\label{examples}

	 \subsection{$SL_2(\CC)$}\label{isom_sl2}
	  	For completeness, we start with the simplest case of 
	  	$G=SL_2(\CC)$. The Cartan subalgebra $\ft$ of diagonal traceless $2\times 2$ matrices is identified with $\CC$ via
	  	$\alpha(A)=A_{11}$ and we
	  	take the $\ZZ/2\ZZ$-invariant polynomial $I=\det$, i.e.,  $I(\alpha)=-\alpha^2$. The cameral (and spectral) curve
	  	$\widetilde{X}_b\subseteq \tot K^2_X$ has equation $\bal^2=\pi^\ast b$, for $b\in \cB=H^0(X,K_X^2)$.
	  	Then, for generic $b$,  the isomorphism (\ref{key_isom})
		\[
			\xymatrix@1{\left(\nabla^{GM}\bla_{SW}\right)_b: &T_b\cB=H^0(X,K_X^2)\ar[r]^-{\simeq} &H^0(\widetilde{X}_b, K_{\widetilde{X}_b})^{\ZZ/2}_{-}  }
		\]

	  	is given by 
	  	\begin{equation}\label{key_isom_sl2}
	  		\bg\longmapsto \left(\nabla^{GM}_\bg\bla_{SW}\right)_b= \left. \frac{\pi^\ast \bg}{2\bal^2}\bla_{SW}\right|_{\widetilde{X}_b} 
	  		= - \left. \frac{\pi^\ast \bg}{2\sigma}\lambda\right|_{\widetilde{X}_b}= 
	  		\left. \frac{\pi^\ast \bg}{2\bla_{SW}}\right|_{\widetilde{X}_b}
	  	\end{equation}
	  	where $\sigma=-\bal^2\in H^0(M,\pi^\ast K_X^2)$ is the tautological section (\ref{taut_section})
	  	of $\bU= K_X^2$   and $\bla_{SW}=\bal$
	  	 is the Liouville (Seiberg--Witten) form.

	  \subsection{$SL_3(\CC)$}
	Consider $G=SL_3(\CC)$, with the standard choices of Borel (upper-triangular) and Cartan (diagonal) subgroups.
	  	Here $\ft\subseteq \fs\fl_3(\CC)$ is the subspace of diagonal traceless
	  	$3\times 3$ matrices
	  	and $W=S_3$.
	  	If we set $\alpha_1(A)=A_{11}-A_{22}$, $\alpha_2(A)=A_{22}-A_{33}$ (two simple positive roots), then we can choose the invariant polynomials to be
	  	\begin{equation}\label{sl3}
	  	\left|
		    \begin{array}{l}
		    I_1(\alpha_1,\alpha_2)= \alpha_1^2 +\alpha_1\alpha_2 +\alpha_2^2\\
		    I_2(\alpha_1,\alpha_2)= -2\alpha_1^3-3\alpha_1^2\alpha_2 +3\alpha_1\alpha_2^2+ 2\alpha_2^3.\\
		    \end{array}
		\right. 
		\end{equation}
	In fact, these are  $I_1(A)= -3(A_{11}A_{22}+A_{11}A_{33}+A_{22}A_{33})$ and $I_2(A)=-27\det A$.
	
	Consequently,  the cameral curve $\widetilde{X}_b$, corresponding to a generic section $b=(b_1,b_2)\in \cB= H^0(X,K_X^2)\oplus H^0(X,K_X^3)$
	is cut out in $M=\tot (K_X^2\oplus K_X^3)$ by the equations
	\begin{equation}\label{sl3_cameral}
		\left|
		    \begin{array}{l}
		    	\bal_1^2+\bal_1\bal_2+\bal_2^2=\pi^\ast b_1\\
		    	-2\bal_1^3-3\bal_1^2\bal_2+3\bal_1\bal_2^2+2\bal_2^3= \pi^\ast b_2\\
		    \end{array}
		\right. 
	\end{equation}

	and
	\[
		\left(\nabla^{GM}_\bg \bla_{SW}\right)_b=\left. 
		\frac{1}{\det \io D\bI}
		\begin{bmatrix}
		       3\bal_1^2- 6\bal_1\bal_2 - 6\bal_2^2 &  2\bal_2 + \bal_1\\
		      -6\bal_1^2 - 6\bal_1\bal_2 + 3\bal_2^2 & - 2\bal_1 -\bal_2 \\
		\end{bmatrix}
		\begin{bmatrix}
			g_1\\ g_2\\
		\end{bmatrix}\right|_{\widetilde{X}_b},
	\]
	    where
	\[
		\det \io D\bI= 27 \bal_1\bal_2 (\bal_1+\bal_2).
	\]

	\subsection{$G_2$}\label{G2_example}
	    It is well-known (\cite{humphreys}[p.103]) that the $G_2$  root system can be embedded in the $B_3$ root system -- and that in fact, 
	    this can be done in a way that simple roots
	of the former are expressed as linear combinations of simple roots of the latter. An explicit description of such an embedding can be obtained by extending
	 the
	calculations in \cite{adlervm_cis}[\S 4], but we do not need  this now.
	      Using this embedding,  we can take  the $\fg_2$ Cartan subalgebra $\ft$ 
	      to consist of  diagonal matrices of the form
		$h=\diag(-a-b,-a,-b,0,b,a,a+b)$,
	for  $a,b\in \CC$. Two simple roots $\alpha_1$, $\alpha_2\in \ft^\vee$ are, e.g.,  $b$ and $a-b$, i.e.,
	\[
		\alpha_1\left(h\right)=  h_{55},\ \alpha_2\left(h\right)= h_{66}-h_{55}.
	\]
	The six positive roots are then $\alpha_1$, $\alpha_2$, $\alpha_1+\alpha_2$, $2\alpha_1+\alpha_2$, $3\alpha_1+\alpha_2$, $3\alpha_1+2\alpha_2$.
	The characteristic polynomial of $h\in \ft$ is
	\[
		\det(h-\lambda E_7)= -\lambda^7+ \lambda^5 2I_1(h) -\lambda^3 I_1^2(h)  + \lambda I_2(h),
	\]
	where, if 
	 we use $\alpha_1$ and $\alpha_2$ as coordinates on  $\ft$,  we have for the invariants
	 \begin{equation}\label{G2_polys}
	 	\left|
		  \begin{array}{l}
		  	I_1(\alpha_1,\alpha_2)= 3\alpha_1^2+3\alpha_1\alpha_2+\alpha_2^2\\
		  	I_2(\alpha_1,\alpha_2)= 4\alpha_1^6+12\alpha_1^5\alpha_2 + 13\alpha_1^4\alpha^2_2+ 6 \alpha_1^3\alpha_2^3+ \alpha_1^2\alpha_2^4.\\
		  \end{array}
	 	\right. 
	 \end{equation}

	The eigenvalues of a matrix from $\fg_2\subseteq \fs\fo_7$ are $0, \pm\lambda_1,\pm \lambda_2,\pm \lambda_3$,
	$\sum_{i=1}^3 \lambda_i=0$. The two invariants are, respectively,  $\frac{1}{2}(\lambda_1^2+\lambda_2^2+\lambda_3^2)$
	and $(\lambda_1\lambda_2\lambda_3)^2$.

	Consequently, the cameral curve $\widetilde{X}_b$, corresponding to a generic section $(b_1,b_2)\in \cB= H^0(X,K_X^2)\oplus H^0(X,K_X^6)$
	is cut out in $M=\tot (K_X^2\oplus K_X^6)$ by the equations

	\begin{equation} \label{G2_cameral}
	 	\left|
		  \begin{array}{l}
		  	3\bal_1^2+3\bal_1\bal_2+\bal_2^2=\pi^\ast b_1\\
		  	4\bal_1^6+12\bal_1^5\bal_2 + 13\bal_1^4\bal^2_2+ 6 \bal_1^3\bal_2^3+ \bal_1^2\bal_2^4= \pi^\ast b_2.\\
		  \end{array}
	 	\right. 
	 \end{equation}
	 
	Now, identifying the adjoint quotient $\chi: \ft\to \ft/W$ with $\bI=(I_1,I_2):\CC^2\to \CC^2$, we obtain  
	that under the isomorphism from Theorem A a section $\bg=(g_1,g_2)^T\in H^0(X,K_X^2)\oplus H^0(X,K_X^6)=T_b\scB$ is mapped to
	$\left(\nabla^{GM}_\bg\bla_{SW}\right)_b$, i.e., 
	\[
		\left. \frac{1}{\det \io D\bI} \begin{bmatrix}
		    	-2\bal_1^2(6\bal_1^3+13\bal_1^2\bal_2+9\bal_1\bal_2^2+2\bal_2^3) &  3\bal_1 + 2\bal_2\\
		    	2\bal_1(12\bal_1^4+30\bal_1^3\bal_2+26\bal_1^2\bal_2^2+9\bal_1\bal_2^3+\bal_2^4) & - 6\bal_1 -3\bal_2
		    	\\
		    \end{bmatrix}
		    \begin{bmatrix}
		    	\pi^\ast g_1\\ \pi^\ast g_2\\
		    \end{bmatrix}\right|_{\widetilde{X}_b},
	\]
	where
	\begin{align}\label{determinant}
		\det \io D\bI= -2\bal_1\bal_2(\bal_1+\bal_2)(2\bal_1+\bal_2)(3\bal_1+\bal_2)(3\bal_1+2\bal_2).
	\end{align}

	We see  that in all examples $\det D\bI$ is a constant multiple of the product of all positive roots. In fact, this follows from
	a classical result of Steinberg \cite{steinberg}. Hence $(\det D\bI)^2$ is proportional to the discriminant $\mathfrak{D}$ of $\fg$ -- the product of all
	roots. Being  $W$-invariant, the discriminant can be expressed as a polynomial in the generators of $\CC[\ft]^W$ -- here,  $I_1$ and $I_2$.
	We sketch a possible way of obtaining this expression without too much brute force. Using the embedding $\fg_2\subseteq \fs\fo_7$,
	we can identify $\alpha_1$ and $\alpha_2$ as $\lambda_1$ and $\lambda_2-\lambda_1$, up to reordering $\lambda_i$'s. Consequently, 
	\begin{equation}
\begin{split}
\mathfrak{D}&= \alpha_1^2(\alpha_1+\alpha_2)^2(2\alpha_1+\alpha_2)^2\alpha_2^2(3\alpha_1+ \alpha_2)^2(3\alpha_1+2\alpha_2)^2\\
&=\lambda_1^2\lambda_2^2\lambda_3^2(\lambda_1-\lambda_2)^2(\lambda_1-\lambda_3)^2(\lambda_2-\lambda_3)^2.
\end{split}
\end{equation}
	The product of the first three terms is $I_2$. The product of the  last three terms is a  polynomial of degree $6$, and 
	hence must be a linear combination of $I_1^3$ and $I_2$, and so one checks
	immediately that
	\[
		\mathfrak{D}=  I_2(4I_1^3-27 I_2).
	\]
	Consequently, the restriction of $\io\mathfrak{D}$ to the universal cameral cover $\cX$
	is the pull-back of a section of $\cO_{\scB}\otimes H^0(X,K_X^{12})$, namely,
	\[
		\scB \ni (b_1,b_2)\longmapsto b_2(4b_1^3-27b_2).
	\]
	It is the Lie derivative of this section that enters the Balduzzi--Pantev formula and its
	generalization (\cite{balduzzi},\cite{ugo_peter_cubic}).
	
	We refer the reader to the beautiful papers \cite{katz_pan_G2},
	\cite{hitchin_G2} for additional details on the $G_2$-Hitchin system, including
	Langlands duality and the description of  Hitchin fibres.

\bibliographystyle{alpha}
\bibliography{biblio}

\end{document}